\documentclass[reqno,12pt]{article}
\usepackage{fullpage,amssymb,amsmath,amsthm,bussproofs}
\usepackage{color,graphics,pdflscape,multicol}
\usepackage[colorlinks]{hyperref}
\usepackage{tikz}
\usetikzlibrary{fit}
\usetikzlibrary{shapes}
\usetikzlibrary{arrows}
\usetikzlibrary{decorations.markings}
\usetikzlibrary{positioning}

\tikzset{mylabel/.style={font=\footnotesize}}
\tikzset{mymidlabel/.style={fill=white}}

\newcommand{\UU}{\mathcal{U}}

\newcommand{\apart}{\,\#\,}
\newcommand{\aut}{\textnormal{\textbf{Aut}}}
\newcommand{\fracfld}{\textnormal{\textbf{Frac}}}
\newcommand{\aparton}{\textnormal{\textbf{Apart}}}
\newcommand{\acommrng}{\textnormal{\textbf{aCRng}}}
\newcommand{\aintdom}{\textnormal{\textbf{aDom}}}
\newcommand{\afld}{\textnormal{\textbf{aFld}}}
\newcommand{\hfiber}{\textnormal{\textbf{hfiber}}}

\newcommand{\iscontr}{\textnormal{\textbf{iscontr}}}
\newcommand{\isweq}{\textnormal{\textbf{isweq}}}
\newcommand{\hset}{\textnormal{\textbf{hSet}}}
\newcommand{\hprop}{\textnormal{\textbf{hProp}}}
\newcommand{\isaprop}{\textnormal{\textbf{isaProp}}}
\newcommand{\isaset}{\textnormal{\textbf{isaSet}}}
\newcommand{\hrel}{\textnormal{\textbf{hRel}}}

\newcommand{\myemph}[1]{\textnormal{\textbf{#1}}}

\newcommand{\id}[1]{\texttt{Id}_{#1}}

\newcommand{\nat}{\mathbb{N}}
\newcommand{\hz}{\mathbb{Z}}
\newcommand{\pathsto}{\leadsto}
\newcommand{\precarry}{\mathbf{p}}
\newcommand{\carry}[1]{#1^{\natural}}
\newcommand{\hexists}{\exists}
\newcommand{\natsummation}{\bigoplus}

\newtheorem{theorem}{Theorem}[section]
\newtheorem{lemma}[theorem]{Lemma}
\newtheorem{proposition}[theorem]{Proposition}

\theoremstyle{definition}
\newtheorem{definition}[theorem]{Definition}
\newtheorem{example}[theorem]{Example}

\theoremstyle{remark}
\newtheorem{remark}[theorem]{Remark}

\title{A preliminary univalent formalization of the $p$-adic numbers}
\author{\'{A}lvaro Pelayo \and Vladimir Voevodsky \and Michael A. Warren}
\date{}

\begin{document}

\maketitle
\setcounter{secnumdepth}{2}   
\setcounter{tocdepth}{1}

\begin{abstract}
The goal of this paper is to give a preliminary formalization
of the $p$-adic numbers, in the context of the Univalent Foundations.  
We also provide the corresponding code verifying the construction
in the proof assistant Coq.  Because work in the univalent setting is ongoing, 
the structure and organization of the construction of the $p$-adic numbers we give in
this paper is expected to change as Coq libraries are more suitably rearranged, and 
optimized, by the authors and other researchers in the future. So our construction here
should be deemed as a first approximation which is subject to improvements.
\end{abstract}

\section{Introduction}

In this paper we present a preliminary formalization of the construction of the
$p$-adic numbers in the Coq proof assistant.  The formalization is
carried out in the \emph{univalent setting} introduced by the second
author \cite{Vo2010}.  This setting, which is based on insights from homotopy theory
and higher-dimensional category theory, serves as an overall
organizational and methodological framework which informs our
construction.  At the same time, our construction has several
ingredients which are familiar in constructive mathematics.
Because work on formalization in this direction is ongoing,
the Coq code associated with this paper may be updated accordingly
in the future by the authors and others.  As such, the structure and
content of the Coq code described here may not match exactly the code
which is ultimately included in the Univalent Foundations libraries.
Readers interested in making use of the code should accordingly
consult the latest version available.

We chose to formalize the $p$-adic numbers as a first step in the
development and formalization of the $p$-adic theory of integrable
systems.  We hope that this will prove to be a promising approach to
this theory which should facilitate progress in the field in the
future, in particular with regard to the construction of algorithms
and their numerical analysis.  Ultimately, we hope that insights from
this project could be useful in the setting of real integrable
systems.

The idea of the univalent perspective is, roughly, to develop
mathematics within the world of homotopy types. By virtue of taking
this approach we are able to make use of type theory as a calculus for
formal reasoning about homotopy types.  We hope that in the future, because this
development of mathematics can be carried out in a proof assistant
such as Coq so that the proofs carry some algorithmic content, it will
be possible to extract good algorithms from the proofs. One of our motivations is that the
construction of such algorithms would in turn help with some problems
concerning integrable systems which are of particular interest in
applications.  For instance, one outstanding problem is: given
numerical spectral data about a quantum system (coming from an experiment), extract an algorithm to
reconstruct the classical integrable system, see Section
\ref{sec:integrable}.

We will only briefly touch upon the technical details of homotopy type
theory and the univalence axiom, and we refer the reader to
\cite{Aw2012} for a basic introduction to homotopy type theory.  For
univalent foundations and the second author's Coq library \cite{Vo2012a} we refer
readers to \cite{PeWa2012}, where a description of the research
program, its motivations, and its implementation in Coq, are given.
Because it is assumed that the reader is already familiar with Coq and
with the second author's program, this paper has been written in a
style which we foresee future papers in formalization taking: it is a
summary of the Coq code written in ordinary mathematical English.  The
details are of course in the Coq code, but the overall structure of
the formalization (as well as the key steps of the proofs) should be
apparent from the sketch given here.  The actual Coq code associated
to this paper can be found on the websites of the authors, as
supplementary files to the arXiv posting of this paper, and as an
appendix to the present paper.

\subsection*{Structure of paper}

Hensel \cite{Hensel} invented the $p$-adic numbers $\mathbb{Q}_{p}$ about
one hundred years ago.  The $p$-adic numbers and the
reals are the canonical metric completions of the rationals.
Classically, there are a number of ways to construct the $p$-adic
numbers, and we refer the reader to \cite{Go1993,Ko1984,
  Sc1984} for further details regarding the classical theory.  The
construction of the $p$-adic numbers given in this paper is
constructive and uses algebraic, rather than analytic, techniques.
Namely, we first construct the integral domain of $p$-adic integers $\hz_{p}$
as a quotient of the ring $\hz[[X]]$ of formal power series over
$\hz$.  We were unable to find the specific construction of $\hz_{p}$
we employ in the literature, but we believe that it is known.  We
then take the $p$-adic numbers $\mathbb{Q}_{p}$ to be the field of
fractions of $\hz_{p}$.  Because we are working constructively, and
because $\hz[[X]]$ does not have decidable equality, it is necessary
to work with an apartness relation and with the corresponding notions
of integral domains and fields.  We will refer to the apartness
versions of fields as \emph{Heyting fields} following the standard
usage in constructive mathematics.

In detail, this paper is organized as follows.  In Section
\ref{sec:basics}, we give a brief overview of the univalent setting.
In Section \ref{sec:constructive_alg} we review some basic
constructive algebra.  Section \ref{sec:fps} contains our construction
of formal power series and the proofs of several
results on formal power series.  The proof that it is
possible to form the Heyting field of fractions for an integral
domain is given in Section \ref{sec:field}.  The construction of the
$p$-adic numbers appears in Section \ref{sec:padics}.  Section
\ref{sec:integrable} is a brief epilogue containing a sketch
of some future plans concerning $p$-adic integrable systems.  Finally, the Coq
code can be found in the Appendix \ref{sec:appendix}.  Although this
appendix is quite long, it is the most important part of the paper and
so we feel that it is justified to include it here.

We should note that the $p$-adic numbers are also relevant in the
physics literature, see \cite{BrFr1993} and the references therein.
In fact, one of our main motivations in wanting to develop a $p$-adic theory of integrable systems
is to study inverse spectral problems concerning $p$-adic analogues
of real quantum integrable systems.  We refer to Section \ref{sec:plans} for a list of short term plans concerning the $p$-adic numbers.

\section{Univalent basics} \label{sec:basics}

The second author's Coq library span a large portion of mathematics
and we make free use of this library.  However, for the sake of
clarity we will here mention those specific parts of the library which
we use in the construction of the $p$-adic numbers.  A survey of the
development of univalent mathematics in Coq can be found in \cite{PeWa2012}.

\subsection*{Notation and conventions}

In this paper, and in the Coq files, all rings are assumed to be
commutative and with $1$.

$\nat$ denotes the type of natural numbers which is defined as an
inductive type in the standard way.  In the Coq code $\nat$ is denoted
by $\mathtt{nat}$.  Similarly, $\hz$ denotes the type of
integers which is constructed as the group completion of the abelian
monoid of natural numbers. In the Coq code $\hz$ is denoted by
$\mathtt{hz}$.

$\UU$ denotes a fixed universe of types.  In the Coq
code this is denoted by $\mathtt{UU}$. The identity type $\id{A}(a,b)$
is denoted by $a\pathsto b$.  In the Coq files this is denoted by
either $\mathtt{paths} \;a \;b$ or by $a\;\sim>\;b$.  

We write $\prod_{x:A}.B(x)$ for dependent
products and $\sum_{x:A}.B(x)$ for dependent sums (defined here as the record type
\texttt{total2}).

We will generally use the same naming conventions as used in the Coq
files, but in some cases we will introduce abbreviations, such as
$\sum_{i=0}^{n}f(i)$ for summation, when it will improve the readability.

Because the current implementation of the underlying type system of
Coq does not handle universes (and several related matters) in a way
which is completely suited for the univalent development of
mathematics, it is necessary to apply several patches to the Coq
system in order to compile the second author's Coq library as well
as the files described in this paper.  Instructions on how to compile
a patched version of Coq can be found in the second author's library.

\subsection{Basic homotopy theoretic notions in Coq}

We think of $\UU$ as the universe of small homotopy types (or fibrant and
cofibrant spaces). For $B:\UU$, we represent a dependent type over $B$
as a term $E:B\to\UU$.  From the perspective of homotopy theory this
corresponds to a fibration over $B$ and, for $b:B$, $E(b)$ corresponds
to the fiber over $b$.  The dependent product $\prod_{x:B}E(x)$ is
regarded as the space of sections of the fibration represented by
$E$.  Similarly, the dependent sum, $\sum_{x:B}E(x)$ corresponds to
the total space of the fibration.  We think of the identity type
$a\pathsto b$ as denoting the fiber of the path space over $(a,b)$.
We will use the phrases "path space" and "type of paths"
interchangeably for this type. I.e., a term $f:a\pathsto b$
corresponds to a path from $a$ to $b$.

Given a path $f:b\pathsto b'$ in $B$ and a point $e:E(b)$ in the fiber
over $b$ we obtain a corresponding point $f_{!}(e):E(b')$ in the fiber
over $b'$.  In the Coq code $f_{!}$ is denoted by $\texttt{transportf}~
E~f~e$.  In order to construct a path $x\pathsto y$ in the total space
$\sum_{x:B}E(x)$ it suffices to construct a path
$f:\pi_{1}(x)\pathsto\pi_{1}(y)$ and a path $g:f_{!}(\pi_{2}(x))\pathsto\pi_{2}(y)$.

Given a term $g:B\to A$ and a path $f:b\pathsto b'$ in $B$, we
obtain a path $g(f):g(b)\pathsto g(b')$.  In the Coq code $g(f)$ is
denoted by $\texttt{maponpaths}~g~f$.  This corresponds, regarding a
homotopy type as an $\infty$-groupoid, the weakly functorial action of $g$ on
the path $f$.

\begin{definition}[\texttt{hfiber}]
  Given types $A$ and $B$, $g:B\to A$ and $a:A$, the \myemph{homotopy fiber of $g$
    over $a$} is the type
  \begin{align*}
    \hfiber\; g\; a  & := \sum_{x:B}\;\bigl(g(x)\pathsto a\bigr).
  \end{align*}
\end{definition}

\begin{definition}[\texttt{iscontr}]
  We define the type $\iscontr(A)$ of proofs that $A$ is contractible
  as 
  \begin{align*}
    \iscontr(A) & := \sum_{c:A}\;\prod_{x:A}\;\bigl(x\pathsto c\bigr).
  \end{align*}
  We say that $A$ is \myemph{contractible} if $\iscontr(A)$ is inhabited.
\end{definition}
We will see below that contractibility in this setting plays the same role as
canonical existence in the classical development of mathematics.
\begin{definition}[\texttt{isweq} and \texttt{weq}]
  Given $g:B\to A$ we define the type $\isweq(g)$ of proofs that $g$
  is a weak equivalence as
  \begin{align*}
    \isweq(g) & := \prod_{x:A}\;\iscontr(\hfiber\;g\;x).
  \end{align*}
  If $\isweq(g)$ is inhabited, then we say that $g$ is a \myemph{weak
    equivalence}.
\end{definition}
There is a filtration of types into different "h-levels".  Homotopy
theoretically this is a slight extension of the usual filtration by
homotopy $n$-types.  We will only require the first few h-levels in
this paper.
\begin{definition}[\texttt{isofhlevel}, \texttt{isaprop},
  \texttt{hprop}, \texttt{isaset} and \texttt{hset}]
  A type $A$ is of \myemph{h-level}:\footnote{Note that in order to
    define $\texttt{isofhlevel}$ as a type which has values in $\UU$,
    as is done in the file \texttt{uu0.v} from the second author's Coq
    library, it is necessary to compile Coq with a patch.}
  \begin{itemize}
  \item $0$ if $A$ is contractible;
  \item $(n+1)$ if, for all $a,b:A$, the type $(a\pathsto b)$ is of h-level $n$.
  \end{itemize}
  We denote by $\iota_{n}(A)$ the type of proofs that $A$ is of
  h-level $n$.  We abbreviate $\iota_{1}(A)$ by $\isaprop(A)$ and
  $\iota_{2}(A)$ by $\isaset(A)$.  We write $\hprop$ for the type of (small) types of
  h-level $1$ and $\hset$ for the type of (small) types of h-level $2$.
\end{definition}
Intuitively, $\hprop$ consists of those spaces which are homotopy
equivalent to either the empty space $0$ or to the one element space
$1$. Accordingly, $\hprop$ plays the role played by the Booleans in
classical logic or by the subobject classifier in topos logic.  Types in $\hprop$ satisfy proof-irrelevance (\texttt{proofirrelevance}) and,
indeed (\texttt{invproofirrelevance}), being an h-prop is equivalent to being proof-irrelevant.

Intuitively, $\hset$ consists of
those spaces which are homotopy equivalent to discrete spaces.  I.e.,
these are the sets.  Most of the types which we will be dealing with
are either h-props or h-sets.  We will sometimes refer to h-sets
simply as "sets" when no confusion will result.

We make use of a number of basic properties of h-levels.  E.g.,
\begin{enumerate}
\item \texttt{impred}: for $n:\nat$, $B:\UU$ and $E:B\to\UU$, the type
  \begin{align*}
    \prod_{x:B}\;\textbf{isofhlevel}_{n}(E_{x})\to\textbf{isofhlevel}_{n}(\prod_{x:B}\;E_{x})
  \end{align*}
  is inhabited.
\item \texttt{impredfun}: for $n:\nat$, $A,B:\UU$,
  if $A$ is of h-level $n$, then so is $(B\to A)$.
\item \texttt{isofhleveldirprod}: If $A$ is of h-level $n$ and $B$ is
  of h-level $n$, then so is $A\times B$.
\end{enumerate}

\subsection{Function extensionality} \label{sec:funext}

We make extensive use of the principle of function extensionality
(\texttt{funextfun}), which follows from the second author's
\emph{Univalence Axiom}.
\begin{definition}[\texttt{funextfun}]
  The principle of \myemph{function extensionality} states that, for
  any two functions $f,g:A\to B$, the type 
  \begin{align*}
    \bigl(\prod_{x:A}\;f(x)\pathsto g(x)\bigr)\to ( f \pathsto g)
  \end{align*}
  is inhabited.
\end{definition}

\subsection{Properties of $\hprop$}

Given a type $A:\UU$, there is a universal way to turn $A$ into a
h-prop.  This is the "inhabited" construction:
\begin{definition}[\texttt{ishinh\underline{~}UU}]
  We say that $A:\UU$ is \myemph{h-inhabited} if the type
  \begin{align*}
    \hat{A} & := \prod_{P:\hprop}((A\to P)\to P)
  \end{align*}
  is inhabited.
\end{definition}
It is immediate, using the facts about h-levels sketched above to see
that $\hat{A}$ is an h-prop.  Moreover, there is a projection
$\pi_{A}:A\to\hat{A}$ given by 
\begin{align*}
  \pi_{A} & := \lambda_{x:A}.\lambda_{P:\hprop}.\lambda_{f:A\to P}.f(x).
\end{align*}
The map $\pi_{A}$ is the universal map from $A$ into a h-prop.  To see
this, observe that if $Q$ is any h-prop and $f:A\to Q$, then we have
a commutative (up to definitional equality) diagram
\begin{align*}
  \begin{tikzpicture}[auto]
    \node (hA) at (0,1.5) {$\hat{A}$};
    \node (A) at (1.5,0) {$A$};
    \node (Q) at (3,1.5) {$Q$};
    \draw[->,bend left=10] (A) to node[mylabel] {$\pi_{A}$} (hA);
    \draw[->,bend right=10] (A) to node[mylabel,swap] {$f$} (Q);
    \draw[->,dashed] (hA) to node[mylabel] {$\bar{f}$} (Q);
  \end{tikzpicture}
\end{align*}
where 
\begin{align*}
  \bar{f} & :=\lambda_{t:\hat{A}}.t(Q)(f).
\end{align*}
Moreover, since $Q$ is a h-prop it follows (using function extensionality) that the space of
such extensions $\bar{f}$ is contractible.

Using the h-inhabited construction it is possible to endow $\hprop$
with the structure of a Heyting algebra.  This structure is summarized
below:
\begin{definition}[\texttt{htrue},\texttt{hfalse},\texttt{hconj},\texttt{hdisj},\texttt{hneg},\texttt{himpl}]
For $P,Q:\hprop$ and $X,Y:\UU$ we define logical operations on
$\hprop$ as follows:
\begin{itemize}
\item $1$ and $0$ are h-props. 
\item $P\wedge Q :=P\times Q$. 
\item $X\vee Y:=\widehat{X+Y}$.
\item $\neg X := X\to 0$. 
\item $X\implies P:=X\to P$.
\end{itemize}
\end{definition}
In addition to the Heyting algebra operations, there is an existential
quantifier (\texttt{hexists}) which is defined by
\begin{align*}
  \hexists_{x:X}P(x) & := \widehat{\sum_{x:X}P(x)}
\end{align*}
for any $P:X\to\UU$ and $X:\UU$.  This quantifier satisfies the usual
properties of the existential quantifier in intuitionistic logic.
Note that our $\hexists$ does \emph{not} correspond to the built-in
existential quantifier "\texttt{exists}" in Coq.

The proof that, with the operations above, $\hprop$ is a Heyting
algebra makes use of the \emph{Propositional
  Univalence Axiom} (\texttt{uahp}) which says that every logical
equivalence between h-props induces a path between them.  I.e., it
says that the type
\begin{align*}
  \prod_{P,Q:\hprop}(P\to Q)\to\bigl( (Q\to P)\to (P\pathsto Q)\bigr).
\end{align*}
is inhabited.

\subsection{Set quotients of types}

The second author has given several constructions of quotients of
types.  A \myemph{hsubtype} of a type $A$ is given by a map
$S\colon A\to\hprop$.  Denote by $\mathcal{P}(A)$ the type of
hsubtypes of $A$.  Given a relation $R$ on
$A$ (that is, $R\colon A\to A\to\hprop$), an \myemph{equivalence class}
consists of a subtype $S$ of $A$ together with the following data:
\begin{enumerate}
\item a term of type $\widehat{\sum_{x:A}S(x)}$.
\item a term of type $\prod_{x,y:A}(xRy \to S(x) \to S(y))$.
\item a term of type $\prod_{x,y:A}(S(x) \to S(y) \to xRy )$.
\end{enumerate}
Given a subtype $S$, we denote by $\textbf{iseqclass}_{R}(S)$ the type
consisting of such data.  The \myemph{set quotient} $A/R$
(\texttt{setquot}) of a type $A$ by a relation $R$ is then defined by
\begin{align*}
  A/R & := \sum_{S:\mathcal{P}(A)}\textbf{iseqclass}_{R}(S).
\end{align*}
It is shown (\texttt{isasetsetquot}) in the second author's library
that $A/R$ is a set and that, when $R$ is an equivalence relation,
this set has the usual universal property.  In particular, there is a
function $\pi\colon A\to A/R$ (\texttt{setquotpr}) which is compatible with the equivalence
relation and, for any set $B$ and function $f:A\to B$ which is
compatible with $R$, there exists an extension $\bar{f}$ making the
diagram
\begin{align*}
  \begin{tikzpicture}[auto]
    \node (hA) at (0,1.5) {$A/R$};
    \node (A) at (1.5,0) {$A$};
    \node (Q) at (3,1.5) {$B$};
    \draw[->,bend left=10] (A) to node[mylabel] {$\pi$} (hA);
    \draw[->,bend right=10] (A) to node[mylabel,swap] {$f$} (Q);
    \draw[->,dashed] (hA) to node[mylabel] {$\bar{f}$} (Q);
  \end{tikzpicture}
\end{align*}
commute.  We will make free use throughout of the results on set
quotients from the second author's library.

\section{Basics on constructive algebra}\label{sec:constructive_alg}

We will here briefly recall some basics of constructive algebra.  For
a more detailed treatment we refer to \cite{BrRi1987} and \cite{MiRiRu1988}.

The usual definitions of fields and integral domains are not entirely
satisfactory from the perspective of constructive algebra since they
deal with negative properties (the property of being a non-zero
element of the field).  From the constructive perspective, it is more
appropriate to replace the notion of an element $x$ being non-zero
($x\neq 0$) with $x$ being {\bf apart from zero}, written $x\apart 0$.  

We will now recall the basics regarding apartness relations.
\begin{definition} (\texttt{isapart}) \label{def:apartness}
  A relation $R:\hrel(X)$ is an \myemph{apartness relation} provided
  that it satisfies the following conditions:
  \begin{description}
  \item[Irreflexive] for all $x:X$, $\neg(xRx)$.
  \item[Symmetric] for all $x,y:X$, $xRy$ implies $yRx$.
  \item[Cotransitive] for all $x,y:X$, if $xRy$, then either $xRz$ or
    $zRy$, for any $z:X$.
  \end{description}
\end{definition}
Classically, the negation of equality $x\neq y$ relation is an
apartness relation.  However, negation of equality is not the only
classical apartness relation.  For example, if $X$ is a topological
space, then the relation $R$ given by $xRy$ if and only if $x$ and $y$
are in different connected components is an apartness relation.  (This
example can be generalized to give a limitless number of classical
examples of apartness relations.)

For $X:\hset$, we denote by $\aparton(X)$ the type of apartness
relations on X. We generally denote apartness relations by $x\apart y$.
When a type has decidable equality the negation of equality is an apartness relation:
\begin{lemma}[\texttt{deceqtoneqapart}] \label{lemm:negeq}
  If $X:\hset$ has decidable equality, then negation of equality
  \begin{align*}
    \neg(x\pathsto y )
  \end{align*}
  is an apartness relation on $X$.
\end{lemma}
\begin{definition}[\texttt{isapartdec}]\label{def:isapartdec}
  Given $X:\hset$ and $R:\aparton(X)$, we say that $R$ is a
  \myemph{decidable apartness relation on $X$} if the type 
  \begin{align*}
    (a R b) + (a \pathsto b)
  \end{align*}
  is inhabited.
\end{definition}
It is immediate (\texttt{isapartdectodeceq}) that if $R$ is a
decidable apartness relation on $X$, then $X$ has decidable equality.

When we are considering algebraic structures equipped with apartness
relations we will require that the relation is compatible with
the operations under consideration.  In particular, for rings we have
the following.
\begin{definition}[\texttt{acommrng}]\label{def:acommrng}
  The type $\acommrng$ consists of commutative rings $A$ together with
  an apartness relation $x\apart y$ on $A$ which is compatible with
  the ring structure of $A$ in the sense that\footnote{Note that in
    the Coq files we actually require the corresponding cancellation
    properties also on the right.  This is redundant for commutative
    rings, but for general rings one requires also these further
    properties.}
  \begin{itemize}
  \item[(i)] 
  For all $a,b,c:A$, if $(c + a )\apart(c+b)$, then $a\apart
    b$.
  \item[(ii)] For all $a,b,c:A$, if $(c\cdot a )\apart(c\cdot b)$, then
    $a\apart b$.
  \end{itemize}
\end{definition}

When a commutative ring $A$ has decidable equality it is
straightforward to verify that negation of equality is compatible with
the ring operations in the sense of Definition \ref{def:acommrng}.

\begin{definition}[\texttt{aintdom}] \label{def:aintdom}
  The type $\aintdom$ consists of $A:\acommrng$ such that 
  \begin{itemize}
  \item $1\apart 0$.
  \item For all $a,b:A$, if $a\apart 0$ and $b\apart 0$, then $(a\cdot
    b)\apart 0$.
  \end{itemize}
\end{definition}
We refer to the terms of type $\aintdom$ as \myemph{apartness
  domains}.

Heyting fields are the appropriate generalization of fields
to the constructive setting when one considers algebraic structures
with apartness relations:
\begin{definition}[\texttt{afld}] \label{def:afld}
  The type $\afld$ of \myemph{Heyting fields} consists of $A:\acommrng$ such that
  \begin{itemize}
  \item $1\apart 0$.
  \item For all $a:A$, if $a\apart 0$, then $a$ has a multiplicative
    inverse (the type of multiplicative inverses of $a$ is inhabited).
  \end{itemize}
\end{definition}
We have the following immediate observation:
\begin{lemma}[\texttt{afldtoaintdom}]
  If $A$ is a Heyting field, then $A$ is an apartness domain.
  \begin{proof}
    It is immediate to prove that, in a Heyting field, if $a$
    has a multiplicative inverse, then it is apart from $0$
    (\texttt{afldinvertibletoazero}).  It follows that $1\apart 0$.
    One can show that if $a$ and $b$ both possess multiplicative
    inverses, then so does their product $a\cdot b$
    (\texttt{multinvmultstable}).  It is then immediate that $(a\cdot
    b)\apart 0$ when $a\apart 0$ and $b\apart 0$.
  \end{proof}
\end{lemma}

\section{Formal power series} \label{sec:fps}

Our treatment of formal power series makes use of function
extensionality, since formal power series over a commutative ring $R$
are here defined as terms of type $\nat\to R$ with the operations of addition
and multiplication given in the usual way. The main result of this
section is that, with these operations, formal power series is
a commutative ring. Moreover, there is a natural apartness relation on
formal power series and, furthermore, when the ring $R$ has decidable
equality the ring of formal power series over $R$ forms an apartness
domain.  We will now fill in the details of this sketch.

\subsection{Summation in a ring} \label{sec:sumation}

We define both a restrictive summation operation
(\texttt{natsummation0}), which allows us to form the sum
$\sum^{n}_{i=0}a_{i}$ of a sequence $a:\nat\to R$, and a more general
operation (\texttt{summation}), which allows us to form the sum
$\sum^{n}_{i=m}a_{i}$ of a sequence $a:\hz\to R$.  However, we will
only really require the former of these two constructions and so we
will omit details related to the more general summation.  In order to
avoid confusion with our notation for dependent sums, we write
$\natsummation_{i=0}^{n}a_{i}$ for the sum $\sum_{i=0}^{n}a_{i}$.
Summation is, of course, defined inductively by setting
\begin{align*}
  \natsummation_{i=0}^{0}a_{i} := a_{0}\quad\text{ and }\quad\natsummation_{i=0}^{n+1}a_{i} := \bigl(\natsummation_{i=0}^{n}a_{i}\bigr)+a_{n+1}.
\end{align*}

\subsubsection*{Manipulation of sums}

It is important to note that when we manipulate sums, to obtain new
sums, \emph{what is relevant is that there is a path between them, and not
whether they are equal in the strict sense}.  This is a crucial point
which underlies in a fundamental way much of the univalent approach to
mathematics.  The following lemma includes several basic facts regarding the behavior of
summation of which we will make frequent use:
\begin{lemma} \label{lemma:key}
  Given a natural number $n$ and sequences $a,b:\nat\to R$, we have
  the following:
  \begin{enumerate}
  \item \textnormal{(\texttt{natsummationpathsupperfixed})} 
    Given $p:\prod_{x:\nat}(x\leq n)\to(a_{x}\pathsto b_{x})$, the type
    \begin{align*}
      \natsummation^{n}_{i=0}a_{i}&\pathsto\natsummation^{n}_{i=0}b_{i}
      \end{align*}
    is inhabited.
  \item \textnormal{(\texttt{natsummationshift0})} 
    The type 
    \begin{align*}
    \natsummation_{i=0}^{n+1}a_{i} &\pathsto (\natsummation^{n}_{i=0}a_{i+1})+a_{0}
    \end{align*}
    is inhabited.
  \end{enumerate}
\end{lemma}
In order to more easily handle reindexing of sums we
introduce, for $f:\nat\to\nat$, the type $\aut_{n}(f)$
(\texttt{isnattruncauto}) of proofs that $f$ is an automorphism of the
interval $[0,n]$ of natural numbers.  Explicitly, $\aut_{n}(f)$ is
defined to be the following type:\footnote{Note that we could,
  alternatively, have used the type $(\prod_{x\leq n}\;\sum_{y\leq
    n}\;(f(y)\pathsto x))\times(\prod_{x\leq n}(f(x)\leq n) )$.
  However, the more verbose type we give here is convenient, for
  purposes of formalization, as it allows for more direct proofs of subsequent lemmas.}
\begin{align*}
  \biggl(\prod_{x\leq n}\;\sum_{y\leq n}\;\bigl(( f(y)\pathsto x )\times
  \prod_{z\leq n}\;(f(z)\pathsto x)\to(y\pathsto z)\bigr)\biggr)\times\bigl(\prod_{x\leq n}\;( f(x)\leq n )\bigr)
\end{align*}
where we have abbreviated $\prod_{x:\nat}\;(x\leq n)\to\cdots$ as
$\prod_{x\leq n}\cdots$ and $\sum_{x:\nat}\;(x\leq
n)\times\cdots$ as $\sum_{x\leq n}\cdots$.  It is possible to
reindex sums along such automorphisms, as shown by the following lemma:
\begin{lemma}
  \textnormal{\texttt{(natsummationreindexing)}} Given a natural
  number $n$ and a map $f:\nat\to\nat$ such that $\aut_{n}(f)$ is
  inhabited, the type
  \begin{align*}
    \natsummation_{i=0}^{n}a_{i} & \pathsto \natsummation_{i=0}^{n}a_{f(i)}
  \end{align*}
  for any sequence $a:\nat\to R$, is inhabited.
\end{lemma}
The final fact regarding summation which we require is the following:
\begin{lemma}
  \label{lemma:natsummationswap}
  \textnormal{\texttt{(natsummationswap)}} Given $f:\nat\to\nat\to R$
  and a natural number $n$, the type
  \begin{align*}
    \natsummation^{n}_{k=0}\natsummation^{k}_{l=0}f(l,k-l) & \pathsto \natsummation^{n}_{k=0}\natsummation^{n-k}_{l=0}f(k,l)
  \end{align*}
  is inhabited.
\end{lemma}

\subsection{The ring of formal power series}

We define, for a type $A$, the type of sequences of elements of $A$
(\texttt{seqson}) as the function space $\nat\to A$.  When $A$ is a set so is $\nat\to
A$ and for $A$ a commutative ring we take $\nat\to A$ as the
underlying set (\texttt{fps}) of the ring of formal power series over
$A$.  If $a$ is a sequence on $A$, then we write $a_{n}:A$ for the
result of evaluating the sequence at the natural number $n$.

\subsubsection*{Ring operations on formal power series}

For a given commutative ring $R$, addition and multiplication of
formal power series are defined as usual by the formulae:
\begin{align*}
  (a+b)_n &:= a_n+b_n\\
(a\cdot b)_{n} & := \natsummation_{k=0}^n a_k b_{n-k}.
\end{align*}
The zero sequence $0$ is given by $0_{n}:=0$ for all natural numbers
$n$ and the sequence $1$ is given by $1_{0}:=1$ and $1_{n+1}:=0$ for
all natural numbers $n$. 
\begin{proposition}[\textnormal{\texttt{fpscommrng}}] \label{theo:main1}
  Let $(R,+,\cdot)$ be a commutative ring. Then the set of sequences
  on $R$ with the operations given above is a commutative ring.
  \begin{proof}
    The proof follows from the facts about summation
    described above.  For example, to prove associativity of
    multiplication, we must show that, for all natural numbers $n$, 
    \begin{align*}
      \natsummation_{i=0}^n \Big(\natsummation_{k=0}^i a_k \cdot b_{i-k}\Big)\cdot c_{n-i} & \pathsto
      \natsummation_{j=0}^n a_j \cdot\Big(\natsummation_{l=0}^{n-j} b_l\cdot c_{(n-j)-l}\Big). 
    \end{align*}
    For this, we reason as follows
    \begin{align*}
      \natsummation_{j=0}^{n}\natsummation_{l=0}^{n-j}a_{j}\cdot(b_{l}\cdot
      c_{(n-j)-l})
      \pathsto
      \natsummation_{l=0}^{n}\natsummation_{j=0}^{l}(a_{l}\cdot b_{k-l})\cdot c_{n-l-(k-l)}
      \pathsto
      \natsummation_{l=0}^{n}\natsummation_{j=0}^{l}a_{l}\cdot(b_{k-l}\cdot c_{n-k}),
    \end{align*}
    where the first path is given by Lemma \ref{lemma:natsummationswap} and
    associativity of multiplication in $R$.  In
    the Coq proof this line of reasoning is put together with generous
    use of Lemma \ref{lemma:key}, (\texttt{funextfun}),
    several minor lemmas such as (\texttt{natsummationtimesdistl}),
    and associativity of $R$ itself.
  \end{proof}
\end{proposition}

\subsection{The apartness relation on formal power series}

Although it is not used in the construction of the $p$-adic numbers,
we mention here some results contained in the Coq files regarding
apartness relations on formal power series.

Assume that $R$ is a commutative ring with an apartness relation.
Then there is an induced apartness relation on $R[[X]]$ given by
setting (\texttt{fpsapart})
\begin{align}\label{eq:fps_apart}
  a\apart b \quad\text{ if and only if }\quad\exists_{n:\nat}.a_{n}\apart b_{n}
\end{align}
for $a,b:R[[X]]$.  This apartness relation is compatible with the ring
operations and so we see that $R[[X]]:\acommrng$
(\texttt{acommrngfps}).

For $R$ an apartness domain, provided that the apartness relation on
$R$ is decidable in the sense of Definition \ref{def:isapartdec}, it is possible to show
that $R[[X]]$ is an apartness domain.
\begin{proposition}[\texttt{apartdectoisaintdomfps}]\label{prop:fps_apart}
  For $R:\aintdom$ with decidable apartness, the commutative ring $R[[X]]$ of formal power
  series is an apartness domain when equipped with the apartness
  relation (\ref{eq:fps_apart}).
\end{proposition}
The proof of Proposition \ref{prop:fps_apart} is a consequence
of the following lemma:
\begin{lemma}[\texttt{leadingcoefficientapartdec}]\label{lemma:leadingcoeff}
  For $R:\aintdom$ and $a:R[[X]]$, if $a_{0}\apart 0$, then for any
  $n:\nat$ and $b:R[[X]]$, if $b_{n}\apart 0$, then $(a\cdot b)\apart 0$.
  \begin{proof}
    The proof is by induction on $n$ and is obvious in the base case.
    The induction case splits into two subcases depending on whether
    $b_{0}\apart 0$ or $b_{0}\pathsto 0$.  In the former case,
    $(a\cdot b)_{0}\apart 0$, whereas in the latter case the claim
    follows by applying the induction hypothesis to the sequence $b':R[[X]]$ given
    by $b'_{n}:=b_{n+1}$.
  \end{proof}
\end{lemma}

\section{The Heyting field of fractions} \label{sec:field}

The construction of the Heyting field of fractions from an apartness
domain is a classical result in constructive algebra due to Heyting
and we therefore give only a brief sketch of the details here.
\begin{definition}[\texttt{aintdomazerosubmonoid}]
  Given $A:\aintdom$, we denote by $\tilde{A}$ the submonoid of $A$ (with
  respect to the multiplicative structure of $A$) consisting of those
  $a:A$ such that $a\apart 0$.
\end{definition}
It follows (\texttt{commrngfrac}) that there exists a commutative ring
$A[\tilde{A}^{-1}]$ obtained by localizing with respect to $\tilde{A}$.
It remains to show there exists an apartness relation on
$A[\tilde{A}^{-1}]$ which makes it into a Heyting field.
\begin{definition}[\texttt{afldfracapartrel0}]
For elements $a,c: A\times\tilde{A}$ we define 
\begin{align*}
  a\apart c \quad\text{ if and only if
  }\quad\bigl((\pi_{1}a)\cdot(\pi_{2}c)\bigr)\apart\bigl((\pi_{1}c)\cdot(\pi_{2}a)\bigr).
\end{align*}
\end{definition}
This relation extends to a relation (\texttt{afldfracapartrel}) on $A[\tilde{A}^{-1}]$ and it is
straightforward to show that it is an apartness relation
(\texttt{afldfracapart}) which is
compatible with the ring structure of $A[\tilde{A}^{-1}]$ (\texttt{afldfrac0}).  For
instance (\texttt{iscotransafldfracapartrelpre}), to see that it is cotransitive
suppose given $(a,a')\apart
(c,c')$ and some $(b,b')$.  Then, by the fact that $A$ is an apartness
domain, we see that $a\cdot c'\cdot b'\apart c\cdot a'\cdot b'$.
Therefore, by cotransitivity of the apartness relation of $A$, we have
that either $a\cdot c'\cdot b' \apart b\cdot a'\cdot c'$ or $b\cdot
a'\cdot c' \apart c\cdot a'\cdot b'$.  In the former case it follows
that $a\cdot b'\apart b\cdot a'$.  I.e., $(a,a')\apart (b,b')$.  In
the latter case it similarly follows that $(b,b')\apart(c,c')$.

Given $a\in A\times\tilde{A}$ such that $a\apart 0$, we have
$\pi_{1}(a)\apart 0$ and therefore, we take $a^{-1}$ to be given by
the pair $(\pi_{2}(a),\pi_{1}(a))$.  This definition extends to a
definition of the inverse of an element apart from $0$ in
$A[\tilde{A}^{-1}]$ and it is straightforward to show that this gives
makes $A[\tilde{A}^{-1}]$ a Heyting field:
\begin{theorem}[\texttt{afldfracisafld}]\label{thm:fld_frac}
  For $A:\aintdom$, with the definitions given above, $A[\tilde{A}^{-1}]$
  forms a Heyting field.
\end{theorem}
We refer to the Heyting field from Theorem \ref{thm:fld_frac} as the
\myemph{Heyting field of fractions of $A$} and we write $\fracfld(A)$ for it.

\section{The $p$-adic numbers} \label{sec:padics}

The $p$-adic numbers were invented about one hundred years ago by German mathematician K. Hensel.

\subsection{Basic number theory}

The following definition is the relation of integer divisibility, 
and is given as a two part definition in the Coq file. The first
part says that, given three integers $n,m,k$, if the product
of $n$ and $k$ is $m$, then $n$ divides $m$. The general
definition starts only with $n$ and $m$, and appeals to
the existence of $k$.

\begin{definition}[$\mathtt{hzdiv0}$ and $\mathtt{hzdiv}$]
Let $n$ and $m$ be integers. We write $n|m$ for the type 
\begin{align*}
  n|m & := \hexists_{k:\hz}.(m \pathsto n\cdot k)
\end{align*}
and we say that $n$ {\bf divides} $m$ when $n|m$ is inhabited.
\end{definition}

The division algorithm is then shown to hold via a series of steps.
First, we prove the division algorithm for natural numbers. Recall
that $\mathtt{pr1}$ and $\mathtt{pr2}$ are defined as projections onto
the base and ``specialization'' to a fiber:

\begin{lemma}[$\mathtt{divalgorithmnonneg}$] \label{lem:danat}
  For $n$ and  $m$ of type $\mathtt{nat}$, with $m$ nonzero, there
  exists a term $qr:(\hz\times\hz)$ such that there is a term of
  type
  \begin{align*}
    n \pathsto\bigl(m\cdot \pi_{1}(qr)\bigr)+\pi_{2}(qr)
  \end{align*}
  and there are proofs that $0\leq \pi_{2}(qr)<m$.
\end{lemma}

The proof of Lemma \ref{lem:danat} is by induction on $n$ with, in the
successor step, a case analysis on whether $(r' + 1 ) < m$ or
$r'\pathsto m$ (that such a case analysis is possible follows from
decidability of equality using \texttt{hzlehchoice} from the second
author's library). The proof of the general division algorithm
is then done by a detailed case analysis (on whether $n$ and $m$ are
negative, non-negative or propositionally equal to $0$):
\begin{theorem}[$\mathtt{divalgorithmexists}$] \label{theo:da}
For $n$ and $m$ of type $\hz$ with $m>0$, the space of terms
$qr:\hz\times\hz$ such that the types $n\pathsto ( m\cdot
\pi_{1}(qr))+\pi_{2}(qr)$ and $0\leq \pi_{2}(qr)<|m|$ are inhabited
is contractible.
\end{theorem}
Here, as throughout, \emph{contractibility} corresponds to
\emph{unique existence} in the traditional setting.  One consequence
of the division algorithm is that we obtain the operations of taking the
quotient and remainder of an integer modulo a non-negative integer
(\texttt{hzquotientmod} and \texttt{hzremaindermod}).  These two
operations will play a role in a number of calculations in the sequel.

In addition to the division algorithm we also obtain the familiar
Euclidean algorithm (again stated in terms of contractibility of an
appropriate space):
\begin{theorem}[$\mathtt{euclideanalgorithm}$] 
  Let $n$ and $m$ be integers with $n \neq 0$. Then the space
  $\mathtt{hzgcd}(n,m)$ of greatest common divisors of $n$ and $m$ is
  contractible.
\end{theorem}
We also obtain a form of the B\'{e}zout lemma:
\begin{lemma}[\textnormal{\texttt{bezoutstrong}}]\label{lemma:bezout}
  For all $m,n:\hz$ such that $n$ is non-zero, the type of
  $ab:\hz\times\hz$ for which there exists a term of type $\gcd(n,m)\pathsto \pi_{1}(ab)\cdot n +
  \pi_{2}(ab)\cdot m$ is inhabited.
\end{lemma}
Given $p:\hz$, the type of proofs that $p$ is a prime is defined by
setting
\begin{align*}
  \texttt{isaprime}(p) := (1<p)\times\bigl((m | p ) \to (m \pathsto 1
  )\vee ( m\pathsto p ) \bigr).
\end{align*}
As a consequence of Lemma \ref{lemma:bezout} we obtain
\begin{theorem}[\textnormal{\texttt{acommrng{\_}hzmod}} and \textnormal{\texttt{ahzmod}}]\label{theorem:zmodp_hfield}
  For non-zero $p$ of type $\hz$, $\hz/p\hz$ is a commutative ring
  with compatible apartness relation.  When $p$ is a prime, $\hz/p\hz$ is a Heyting field.
\end{theorem}
Note that the apartness relation on $\hz/p\hz$ is the one induced by the
fact that equality of $\hz/p\hz$ is decidable (\texttt{isdeceqhzmodp}).

\subsection{The construction of $\mathbb{Q}_{p}$}

Throughout this section we assume given a prime $p$.  Explicitly, we
require the proof witnessing the fact that $p$ is a prime.
We note though that for some of the results stated here it is only
necessary that $p$ be non-zero.  We also introduce some
notation for quotients and remainders modulo $p$.  We denote by
$\{a\}$ the quotient of $a$ modulo $p$ (\texttt{hzquotientmod}) and
by $[ a ]$ the remainder of $a$ modulo $p$.

We will now summarize our construction of the apartness domain $\hz_{p}$ of p-adic integers.
\begin{definition}[\textnormal{\texttt{precarry}}]
  Given a formal power series $a$ over $\hz$, we define a new
  formal power series $\precarry(a)$ over $\hz$ inductively by 
  \begin{align*}
    \precarry(a)_{0} & := a_{0}\\
    \precarry(a)_{n+1} & := a_{n+1}+\bigl\{\precarry(a)_{n}\bigr\}.
  \end{align*}
\end{definition}
\begin{definition}[\textnormal{\texttt{carry}}]
  Given a formal power series $a$ over $\hz$, we define a new formal
  power series $\carry{a}$ over $\hz$ by
  \begin{align*}
    (\carry{a})_{n} & := \bigl[\precarry(a)_{n}\bigr].
  \end{align*}
  We call $\carry{a}$ the \myemph{carried power series of $a$}.
\end{definition}
\begin{example}
  The formal power series $a=(4,1,8,0,\ldots)$ is sent to
  $\precarry(a)=(4,2,8,2,0,\ldots)$ and to $\carry{a}=(1,2,2,2,0,\ldots)$.
\end{example}
The operation of carrying (mod $p$) for power series induces an
equivalence relation $\sim$ (\texttt{carryequiv}) on $\hz[[X]]$ by setting 
\begin{align*}
  a \sim b \quad\textnormal{ if and only if }\quad\carry{a}\pathsto \carry{b}.
\end{align*}
Observe that $X-p\sim 0$.  Furthermore, for any $a\in\hz[[X]]$, if
$a\sim 0$, then there exist integers $\lambda_{i}$ such that
$a_{0}=-\lambda_{0}p$ and $a_{n+1}=-\lambda_{n+1}p+\lambda_{n}$.
Using these facts it follows that $\sim$ is the equivalence relation
corresponding to the ideal $(X-p)$ in $\hz[[X]]$.  Ultimately,
once the theory of ideals has been developed in the Univalent
Foundations Library, $\hz_{p}$ will be constructed as the quotient of
$\hz[[X]]$ by this ideal.  However, because quotients of rings are
given in the second author's library in terms of congruences, we here describe
$\hz_{p}$ using the corresponding congruence $\sim$.

We will now describe the proof that this relation is a congruence with
respect to the ring operations on $\hz[[X]]$.
\begin{lemma}[\texttt{quotientprecarryplus}]\label{lemma:precarryandplus}
  For formal power series $a$ and $b$ over $\hz$, 
  \begin{align*}
    \bigl\{\precarry(a+b)_{n}\bigr\} & \pathsto \bigl\{\precarry(a)_{n}\bigr\}+\bigl\{\precarry(b)_{n}\bigr\}+\bigl\{\precarry(\carry{a}+\carry{b})_{n}\bigr\}
  \end{align*}
  for $n:\nat$.
  \begin{proof}
    The proof is by induction on $n$.  In the base case it is
    trivial and in the induction case it is by the following
    calculation:
    \begin{align*}
      \bigl\{\precarry(a+b)_{n+1}\bigr\} & \pathsto
      \bigl\{\precarry(a)_{n+1}+\precarry(b)_{n+1}+\{\precarry(a^{\natural}+b^{\natural})_{n}\}\bigr\}\\
      & \pathsto
      \bigl\{\precarry(a)_{n+1}\bigr\}+\bigl\{\precarry(b)_{n+1}\bigr\}+\bigl\{\precarry(a^{\natural}+b^{\natural})_{n}\bigr\}+\bigl\{a^{\natural}_{n+1}+b^{\natural}_{n+1}+[\precarry(a^{\natural}+b^{\natural})_{n}]\bigr\}\\
      & \pathsto\bigl\{\precarry(a)_{n+1}\bigr\}+\bigl\{\precarry(b)_{n+1}\bigr\}+\bigl\{\precarry(a^{\natural}+b^{\natural})_{n+1}\bigr\}
    \end{align*}
    where the first path is by definition of precarry and the
    induction hypothesis, the second path is by the familiar
    decomposition of the quotient of a sum, and the final path is by
    definition and the fact that the quotient of a remainder is zero.
  \end{proof}
\end{lemma}
The following observation is a consequence of
Lemma \ref{lemma:precarryandplus}.
\begin{lemma}[\texttt{carryandplus}]\label{lemma:carryandplus}
  For $a$ and $b$ formal power series over $\hz$,
  $\carry{(a+b)}\pathsto\carry{(\carry{a}+\carry{b})}$.
\end{lemma}
Similarly, a straightforward induction gives us the following lemma:
\begin{lemma}[\texttt{precarryandtimesl}]
  \label{lemma:precarryandtimes}
  Given formal power series $a$ and $b$ over $\hz$,
  \begin{align*}
    \bigl\{\precarry(a\cdot b )_{n}\bigr\} & \pathsto
    (\{\precarry(a)\}\cdot b )_{n}+\bigl\{\precarry(a^{\natural}\cdot b)_{n}\bigr\}
  \end{align*}
  for $n:\nat$.
\end{lemma}
The proof that carrying is compatible with multiplication of
power series is then an immediate consequence of Lemma \ref{lemma:precarryandtimes}:
\begin{lemma}[\texttt{carryandtimes}]\label{lemma:carryandtimes}
  Given formal power series $a$ and $b$ over $\hz$, $(a\cdot
  b)^{\natural}\pathsto ( a^{\natural}\cdot b^{\natural})^{\natural}$.
\end{lemma}
It follows from Lemmas \ref{lemma:carryandplus} and
\ref{lemma:carryandtimes} that the quotient of $\hz[[X]]$ by the
equivalence relation $\sim$ is itself a commutative ring
(\texttt{commrngofpadicints}).  Indeed, it is the commutative ring
$\hz_{p}$ of \myemph{$p$-adic integers}.  Moreover, there is an apartness
relation (\texttt{padicapart}) on $p$-adic integers obtained as the extension of the
relation (\texttt{padicapart0})
\begin{align}\label{eq:padicapart}
  a \apart b \quad \text{ if and only if }\quad \exists_{n:\nat}.\neg (
  a^{\natural}_{n}\pathsto b^{\natural}_{n}),
\end{align}
for $a,b:\hz[[X]]$, to the $p$-adic integers.  This apartness relation
is straightforwardly seen to be compatible with the ring structure of
$\hz_{p}$ (\texttt{acommrngofpadicints}).
\begin{theorem}[\texttt{padicintsareintdom},\texttt{padicintegers}]\label{thm:padic_ints}
  The commutative ring $\hz_{p}$ with the apartness relation described
  above forms an apartness domain.
  \begin{proof}
    It suffices to prove that for $a,b:\hz[[X]]$ such that $a\apart 0$
    and $b\apart 0$ it follows that $a\cdot b\apart 0$, where we are
    considering only the apartness relation (\ref{eq:padicapart}).
    Since $\hz$ has decidable equality, it follows
    (\texttt{leastelementprinciple}) that there are natural numbers
    $k$ and $m$ which are the least natural numbers such that $\neg
    (a^{\natural}_{k}\pathsto 0)$ and $\neg(b^{\natural}_{m}\pathsto
    0)$, respectively.  It then follows that $\neg((a\cdot
    b)^{\natural}_{k+m}\pathsto 0)$.  

    To see this, assume for a contradiction that there is a path $(a\cdot
    b)^{\natural}_{k+m}\pathsto 0$ and consider first the
    case where $k+m=0$. Then we have that $a_{0}\cdot b_{0}$ is
    congruent to $0$ modulo $p$ and therefore, since $p$ is prime,
    either $a_{0}$ is congruent to $0$ modulo $p$ or $b_{0}$ is
    congruent to $0$ modulo $p$.  In either case we have obtained a
    contradiction.  

    On the other hand, when $k+m$ is a successor $k+m=n+1$, we have that
    \begin{align}\label{eq:padicints}
      (a\cdot b)^{\natural}_{k+m} & \pathsto \bigl[(a^{\natural}\cdot
      b^{\natural})_{k+m}+\{\precarry(a^{\natural}\cdot b^{\natural})^{\natural}_{n}\}\bigr].
    \end{align}
    By the choice of $k$ and $m$ it follows that there is a further
    term (\texttt{precarryandzeromult}) of type
    $\precarry(a^{\natural}\cdot b^{\natural})^{\natural}_{n}\pathsto
    0$.  Therefore, we obtain a term of type
    \begin{align*}
      0 & \pathsto [(a^{\natural}\cdot b^{\natural})_{k+m}].
    \end{align*}
    However, it is easy (\texttt{hzfpstimeswhenzero}) to see that
    $(a^{\natural}\cdot b^{\natural})_{k+m}\pathsto
    a^{\natural}_{k}\cdot b^{\natural}_{m}$.  So, since $p$ is prime,
    either $a^{\natural}_{k}\pathsto 0$ or $b^{\natural}_{m}\pathsto
    0$ is inhabited.  In either case we obtain a contradiction. 
  \end{proof}
\end{theorem}
Using Theorem \ref{thm:padic_ints}, we now arrive at our definition of
the $p$-adic numbers:
\begin{definition}[\texttt{padics}]
  The Heyting field $\mathbb{Q}_{p}$ of $p$-adic numbers is defined as
  the Heyting field of fractions of $\mathbb{Z}_{p}$:
  \begin{align*}
    \mathbb{Q}_{p} & := \fracfld(\mathbb{Z}_{p}).
  \end{align*}
\end{definition}

\section{Future directions: towards $p$-adic integrable systems}
\label{sec:integrable}

Next we present an outline of the work on $p$-adic integrable systems that we plan to
carry out following this paper. The long term goal
is to develop an analogue of the symplectic theory of
finite-dimensional real integrable systems in \cite{PeVu2009, PeVu2011}  for $p$-adic integrable systems in the univalent setting, and implement it in Coq. 

We are beginning to explore this, and what we give next is a brief and informal 
glimpse of our plans. At this point this section is a discussion without rigorous descriptions
as we are not yet convinced of the optimal definition of $p$-adic integrable
system.  We hope to convey the fact that the $p$-adic and
real theories are expected to be different, and draw attention to the
topic; in fact, we are not aware of a uniform treatment of $p$-adic integrable
systems in the symplectic setting.

\subsection{Definition of $p$-adic integrable systems}

\subsubsection*{A word on the contrast between $p$-adic and real notions}

We refer to \cite[Section~3]{Go1993} for basic algebraic and topological aspects
concerning the $p$-adic numbers. Many aspects do not match the intuition we
have for the real numbers. For instance, there are no nontrivial connected
sets and there are non-empty sets which are both compact and open. Other
aspects are more familiar: on $\mathbb{Q}_p$ there is an absolute value
$| \cdot |$ and $\mathbb{Q}_p$ is complete with respect to it, and there is
an inclusion $\mathbb{Q} \to \mathbb{Q}_p$ with dense image. Continuity and
differentiability of functions is defined in the usual way \cite[Definitions 4.2.1, 4.2.2]{Go1993}.
Continuous functions are uniformly continuous on compact sets, as in the real case. 

The notions of continuity and differentiability extend to functions $f \colon U \subset (\mathbb{Q}_p)^n \to \mathbb{Q}_p$
of several variables $(x_1,\ldots,x_n)$ on open sets $U$ of the
Cartesian product $(\mathbb{Q}_p)^n$, in direct analogy with the real case, and in
particular we have analogous definitions for partial derivatives
$\frac{\partial f}{\partial x_i},$ for all $i=1,\ldots,n$. 
But although
the definitions are the same, differentiability behaves differently in the $p$-adic
case than in the real case. For instance, there are functions $f \colon \mathbb{Q}_p \to \mathbb{Q}_p$
which have zero derivative everywhere but are \emph{not} locally constant. Also, the
natural extension of the real mean value theorem to the $p$-adic case is false in
general (although a version holds for sufficiently close points), as seen for instance
by considering $f(x)=x^p-x$ between the extreme points $a=0$ and $b=1$. In this case, \cite[Proposition 4.2.3]{Go1993}
$f'(x)=px^{p-1}-1$ and $f(a)=f(b)=0$ and it is easy to check that any element ``in between''
$a$ and $b$, that is, of the form $at +b(1-t)=1-t$ for some $t$ with $|t|\leq 1$, gives rise
to a unit $f'(1-t)$ in $\mathbb{Z}_p$.

These differences are an indication that the theory of $p$-adic integrable systems 
is not expected to be a direct extension of the theory of real integrable systems, even if the
basic definitions are analogous. One can explore such theory classically only, but we hope to
do it in the univalent setting, building on the constructions of $\mathbb{Q}_p$ which we
have given in the previous sections.

\subsubsection*{Integrable systems}

We are here going to propose a notion of $p$-adic integrable systems
in parallel with the commonly accepted notion of real integrable
systems, at least in symplectic geometry.

Because in the univalent foundations, and in Coq, it is nontrivial to define manifolds, for now we are going to work with the $p$-adic Cartesian
product
\begin{align*}
M:=(\mathbb{Q}_p)^{2n}=\mathbb{Q}_p \times \ldots  
(2n \,\, \textup{times}) \ldots \times \mathbb{Q}_p
\end{align*}
with coordinates $(x_1,y_1,\ldots,x_n,y_n)$.  In this way, we also avoid a discussion of differential or symplectic forms.
Fix a $p$-adic measure on $\mathbb{Q}_p$, and endow $M$ with the induced product 
measure.  

On $M$ we may consider differentiable functions in the $p$-adic sense\footnote{for now 
we are thinking only of polynomials on $2n$-variables, which are easy to deal with in Coq.}.   The
following is the formal extension of the definition of real integrable system in finite dimensions.
There is, however, a critical point which is not clear to us at the moment, and that's
why we restrict our definition to analytic maps, see Remark \ref{analytic}.

\begin{definition} \label{def:psystem}
We will say that a ($p$-adic) analytic map 
$F:=(f_1,\ldots,f_n) \colon M \to (\mathbb{Q}_p)^n$
is a {\bf $p$-adic integrable system} if two conditions hold: 
\begin{enumerate}
\item
the collection $f_1,\ldots,f_n$ satisfies Hamilton's equations:
\begin{align} \label{eq:t}
 \sum_{k=1}^n \frac{\partial f_i}{\partial x_k} \,\frac{\partial f_j}{\partial y_k}
-  \frac{\partial f_i}{\partial y_k} \,\frac{\partial f_j}{\partial x_k} =0, \,\,\,\,\,\,\, \forall \,\,1 \leq i \le j \leq n.
\end{align}
\item
the  set where the $n$ formal differentials 
\begin{align*}
{\rm d}p_i & := \Big(\frac{\partial f_i}{\partial x_1},\ldots,\frac{\partial f_i}{\partial x_n},
\frac{\partial f_i}{\partial y_1},\ldots,\frac{\partial f_i}{\partial y_n}     \Big),\,\,\,\,\,\, \forall \,\,1\leq i \leq n 
\end{align*}
are linearly dependent  has $p$-adic measure $0$. 

That is, there exists a $p$-adic measure $0$ set $A$ such that 
 ${\rm d}f_1,\ldots,{\rm d}f_n$
are linearly independent on $M \setminus A$.  The points where 
 ${\rm d}f_1,\ldots,{\rm d}f_n$ are linearly dependent are 
 called \emph{singularities}.
\end{enumerate}
\end{definition}

\begin{remark} \label{analytic}
  This remark explains why we have to restrict to analytic functions in
  Definition \ref{def:psystem}, when in the real theory one likes
  to include all smooth functions in the definition of integrable system.
  There are many interesting, nontrivial $p$-adic functions that are smooth 
  and have zero derivative everywhere. \emph{However this is not possible if one
    restricts to analytic functions}.   
  Therefore if $f$ is a smooth solution to a linear differential equation, we
  could add to $f$ any of these nontrivial functions with zero derivative and 
  obtain a new solution.   It follows that all collections
  of $n$ smooth functions $f_1,\ldots,f_n$ which are smooth and have zero
  derivative everywhere would also form a kind of integrable system,
  but a very "degenerate" one (in the sense that the differentials
  ${\rm d}f_1, \ldots, {\rm d}f_n$ would not be linearly independent almost everywhere
  as it is normally required for real integrable systems). So this
  undesirable case does not occur.
  However, adding functions with zero derivative to an existing system
  would be unavoidable, giving rise to a  new, seemingly very different,
  $p$-adic integrable system.  We currently understand neither what
  this means geometrically, nor what it implies for the development of
  the theory.
\end{remark}

\subsection{Future plans} \label{sec:plans}

The following is a rough outline of what we would like to do next. 

\subsubsection*{Towards $p$-adic symplectic geometry}

\begin{itemize}
\item[$\rhd$]
\emph{$p$-adic manifolds}: formalize the notion of $p$-adic manifold in the univalent Foundations 
with Coq.  Formalize Serre's theorem \cite{Se1965} classifying compact $p$-adic
manifolds. 

\item[$\rhd$]
\emph{$p$-adic symplectic forms}: a $p$-adic symplectic form $\omega$ may be defined as in the
real case. The closedness condition ${\rm d}\omega=0$ makes sense in the $p$-adic setting,
and so does the non-degeneracy condition (in fact, over any field).  In the real setting,
a theorem of Darboux says that all symplectic forms are locally equivalent, so
real symplectic manifolds have no local invariants.  It is natural to wonder whether
this result holds in the $p$-adic setting "as is". Because of our previous discussion 
(see Remark \ref{analytic}) one should probably restrict to the analytic setting since
${\rm d}\omega=0$ is in fact a system of partial differential equations. Darboux's theorem
plays a leading role in the theory of real integrable systems.
\end{itemize}

\subsubsection*{Towards $p$-adic integrable systems: basic theory}

\begin{itemize}
\item[$\rhd$]
\emph{ construction of $p$-adic integrable systems}: define $p$-adic integrable systems on $p$-adic manifolds, not just 
$(\mathbb{Q}_p)^n$, and implement this
in the univalent foundations using Coq.
\item[$\rhd$]
\emph{ $p$-adic local and semiglobal theory}: 
develop the local and semilocal theory of $p$-adic integrable systems in Coq. The
local theory basically refers to local models, and the semilocal theory refers to local
models in neighborhoods of fibers. One is interested in both the topological
and symplectic classification of such models.   We are not aware of results describing the 
topological, or symplectic, structure  of regular or singular fibers
in the $p$-adic setting.

In the real case, the
regular fibers and their neighborhoods are understood (this is the famous 
Action-Angle Theorem due to Mineur and Arnold.) The singular fibers
may be complicated and not are yet well understood in the real setting either
(if one restrict to the real analytic setting, then the theory is better understood).
\end{itemize}

\subsubsection*{Towards $p$-adic toric and semitoric systems}

\begin{itemize}

\item[$\rhd$]
\emph{ $p$-adic toric systems}: 
a particular class of real integrable systems which has been thoroughly studied and is
well understood, is that of toric integrable systems $F=(f_1,\ldots,f_n)$ on
$2n$-dimensional compact symplectic manifolds $(M,\omega)$. 
These are systems  in which each component $f_i$ generates a flow which is
periodic of a fixed period. In this case, $F$ is called a \emph{momentum map}. 
Atiyah \cite{At1982}, Guillemin-Sternberg \cite{GuSt1982} 
and Delzant \cite{De1988} proved
a series of striking theorems concerning these systems in the 1980s,
which in particular led to complete combinatorial classification in terms of convex
polytopes by Delzant (these convex polytopes are nothing by the images of $M$
under $F$).  A theorem of Serre \cite{Se1965} classifies compact $p$-adic
analytic varieties. If on these varieties we would consider actions of the
$p$-adic $n$-torus, we do not know to what extent the above
results could be extended. If in
Definition \ref{def:psystem} one allows smooth non-analytic functions, these
results would not hold (see Remark \ref{analytic}).

\item[$\rhd$]
\emph{ $p$-adic semitoric systems}: give a classification of $p$-adic integrable systems 
under some periodicity condition in analogy with \cite{PeVu2009, PeVu2011}.  
\end{itemize}

\subsubsection*{Spectral questions for $p$-adic integrable systems}

Here we restrict to the systems in the previous section, for which we know that
in the real case a full classification may be given.

\begin{itemize}
\item[$\rhd$]
\emph{ Inverse spectral problems}:  
construct algorithms to solve inverse spectral problems about quantum integrable
systems. The leading question in the real case is: given the spectrum, can one recover the system from it? 

\item[$\rhd$]
\emph{ Numerical implementation of inverse spectral problems}:  
constructing numerically accurate algorithms to solve inverse spectral problems.
\end{itemize}

The first subsection above should be  within reach. We expect the second and third subsections to
be substantial. The fourth one depends on the third and it is
difficult to predict how complicated it will be.   

\subsubsection*{Acknowledgements.} 
We thank Mark Goresky, Helmut Hofer, Gopal Prasad and Bas Spitters for discussions.

AP was partly supported by NSF Grant DMS-0635607, an NSF CAREER Award
DMS-1055897, Spain Ministry of Science Grant MTM 2010-21186-C02-01, and 
Spain Ministry of Science Sev-2011-0087.  MAW is supported by the
Oswald Veblen fund and also received support from NSF
Grant DMS-0635607 during the preparation of this paper.  VV is supported by NSF Grant DMS-1100938.  This
material is based upon work supported by the National Science
Foundation.  Any opinions, findings and conclusions or recommendations
expressed in this material are those of the authors and do not
necessarily reflect the views of the National Science Foundation.

\section*{Appendix: Coq code} \label{sec:appendix}

\emph{Disclaimer:} The libraries summarized and listed below are
in preliminary form and are actively being improved and extended by
the authors and others.  As such, we advise interested readers to
consult also with the most recent versions, which need not agree in
form and content with the libraries described here.

\begin{figure}[h]
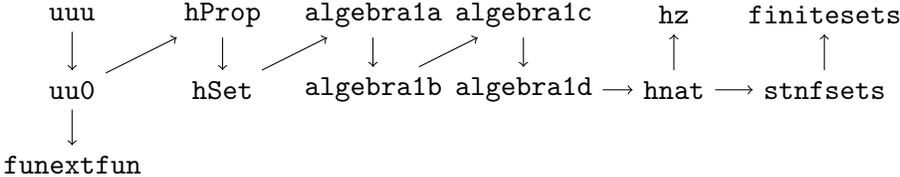

  \centering

  \caption{Dependence diagram of the second author's Coq library}
  \label{figure:VVlib}
\end{figure}

\begin{figure}[h]
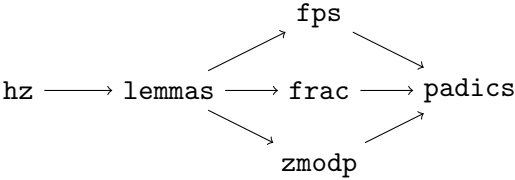

  \centering
      %
  \caption{Dependence diagram of the additional Coq files for the $p$-adics}
  \label{figure:newlib}
\end{figure}

The Coq code is included in full below for easy reference by the
reader.  We also expect to make it available on the webpages of the
first and third authors.  For easy reference, we include here a brief
sketch of the contents of each of the files.  \emph{It is worth remarking
that all of the files described here rely upon the second author's
Coq library.}  For more on this library we refer the reader to
the library itself and to the tutorial \cite{PeWa2012}.  For
quick reference, Figures \ref{figure:VVlib} and \ref{figure:newlib}
give the dependences of the second author's library and the
library associated with this paper, respectively.  

Of the new files, the file \texttt{lemmas.v} contains a number of
small lemmas which, such as basic facts about apartness relations,
some lemmas on rings, \emph{et cetera}, which are required by the
other files.  The file \texttt{fps.v} contains all of the material on
formal power series.  The construction of the Heyting field of
fractions can be found in \texttt{frac.v}.  The basic number theoretic
results which we require are in \texttt{zmodp.v}.  Finally, the
construction of the $p$-adic numbers is given in \texttt{padics.v}. 

\newpage

\begin{landscape}

\subsection{The file \texttt{lemmas.v}}

\begin{multicols}{2}{
{
\tiny
%
}
}\end{multicols}

\subsection{The file \texttt{fps.v}}

\begin{multicols}{2}{
{
\tiny
%
}
}\end{multicols}

\subsection{The file \texttt{frac.v}}

\begin{multicols}{2}{
{
\tiny
\begin{verbatim}
(** *The Heyting field of fractions for an apartness domain *)

(** By Alvaro Pelayo, Vladimir Voevodsky and Michael A. Warren *)

(** February 2011 and August 2012 *)

(** Settings *)

Add Rec LoadPath "../Generalities".  Add Rec LoadPath "../hlevel1".
Add Rec LoadPath "../hlevel2".  Add Rec LoadPath "../Algebra".

Unset Automatic Introduction. (** This line has to be removed for the
file to compile with Coq8.2 *)

(** Imports *)

Require Export lemmas.

(** * I. The field of fractions for an integrable domain with an
apartness relation *)

Open Scope rng_scope.

Section aint.

Variable A : aintdom.

Ltac permute := solve [ repeat rewrite rngassoc2; match goal with | [
  |- ?X ~> ?X ] => apply idpath | [ |- ?X * ?Y ~> ?X * ?Z ] => apply
  maponpaths; permute | [ |- ?Y * ?X ~> ?Z * ?X ] => apply (
  maponpaths ( fun x => x * X ) ); permute | [ |- ?X * ?Y ~> ?Y * ?X ]
  => apply rngcomm2 | [ |- ?X * ?Y ~> ?K ] => solve [ repeat rewrite
  <- rngassoc2; match goal with | [ |- ?H ~> ?V * X ] => rewrite (
  @rngcomm2 A V X ); repeat rewrite rngassoc2; apply maponpaths;
  permute end | repeat rewrite rngassoc2; match goal with | [ |- ?H ~>
  ?Z * ?V ] => repeat rewrite <- rngassoc2; match goal with | [ |- ?W
  * Z ~> ?L ] => rewrite ( @rngcomm2 A W Z ); repeat rewrite
  rngassoc2; apply maponpaths; permute end end ] |[ |- ?X * ( ?Y * ?Z
  ) ~> ?K ] => rewrite ( @rngcomm2 A Y Z ); permute end | repeat
  rewrite <- rngassoc2; match goal with | [ |- ?X * ?Y ~> ?X * ?Z ] =>
  apply maponpaths; permute | [ |- ?Y * ?X ~> ?Z * ?X ] => apply (
  maponpaths ( fun x => x * X ) ); permute | [ |- ?X * ?Y ~> ?Y * ?X ]
  => apply rngcomm2 end | apply idpath | idtac "The tactic permute
  does not apply to the current goal!" ].

Lemma azerorelcomp ( cd : dirprod A ( aintdomazerosubmonoid A ) ) ( ef
: dirprod A ( aintdomazerosubmonoid A ) ) ( p : ( pr1 cd ) * ( pr1 (
pr2 ef ) ) ~> ( ( pr1 ef ) * ( pr1 ( pr2 cd ) ) ) ) ( q : ( pr1 cd ) #
0 ) : ( pr1 ef ) # 0.  Proof.  intros. change ( ( @op2 A ( pr1 cd ) (
pr1 ( pr2 ef ) ) ) ~> ( @op2 A ( pr1 ef ) ( pr1 ( pr2 cd ) ) ) ) in p.
assert ( ( @op2 A ( pr1 cd ) ( pr1 ( pr2 ef ) ) ) # 0 ) as v. apply
A. assumption. apply ( pr2 ( pr2 ef ) ). rewrite p in v. apply ( pr1 (
timesazero v ) ).  Defined.
 
Lemma azerolmultcomp { a b c : A } ( p : a # 0 ) ( q : b # c ) : a * b
# a * c.  Proof.  intros. apply aminuszeroa. rewrite <-
rngminusdistr. apply ( pr2 A ). assumption.  apply
aaminuszero. assumption.  Defined.

Lemma azerormultcomp { a b c : A } ( p : a # 0 ) ( q : b # c ) : b * a
# c * a.  Proof.  intros. rewrite ( @rngcomm2 A b ). rewrite (
@rngcomm2 A c ). apply ( azerolmultcomp p q ).  Defined.

Definition afldfracapartrelpre : hrel ( dirprod A (
aintdomazerosubmonoid A ) ) := fun ab cd : _ => ( ( pr1 ab ) * ( pr1 (
pr2 cd ) ) ) # ( ( pr1 cd ) * ( pr1 ( pr2 ab ) ) ).

Lemma afldfracapartiscomprel : iscomprelrel ( eqrelcommrngfrac A (
aintdomazerosubmonoid A ) ) ( afldfracapartrelpre ).  Proof.  intros
ab cd ef gh p q. unfold afldfracapartrelpre.  destruct ab as [ a b
]. destruct b as [ b b' ].  destruct cd as [ c d ]. destruct d as [ d
d' ].  destruct ef as [ e f ]. destruct f as [ f f' ].  destruct gh as
[ g h ]. destruct h as [ h h' ].  simpl in *.
  
  apply uahp. intro u. apply p. intro p'. apply q. intro q'.  destruct
  p' as [ p' j ]. destruct p' as [ i p' ]. destruct q' as [ q' j'
  ]. destruct q' as [ i' q' ].  simpl in *.
  
  assert ( a * f * d * i * h * i' # e * b * d * i * h * i') as v0.
  assert ( a * f * d # e * b * d ) as v0. apply azerormultcomp. apply
  d'.  assumption.  assert ( a * f * d * i # e * b * d * i ) as
  v1. apply azerormultcomp.  assumption. assumption.  assert ( a * f *
  d * i * h # e * b * d * i * h ) as v2. apply azerormultcomp.  apply
  h'. assumption.  apply azerormultcomp. assumption. assumption.
  apply ( pr2 ( acommrng_amult A ) b ). apply ( pr2 ( acommrng_amult A
  ) f ). apply ( pr2 ( acommrng_amult A ) i ). apply ( pr2 (
  acommrng_amult A ) i' ).

  assert ( a * f * d * i * h * i' ~> c * h * b * f * i * i' ) as l.
  assert ( a * f * d * i * h * i' ~> a * d * i * f * h * i' ) as l0.
  change ( @op2 A ( @op2 A ( @op2 A ( @op2 A ( @op2 A a f ) d ) i ) h
  ) i' ~> @op2 A ( @op2 A ( @op2 A ( @op2 A ( @op2 A a d ) i ) f ) h )
  i' ).  permute. rewrite l0. rewrite j.  change ( @op2 A ( @op2 A (
  @op2 A ( @op2 A ( @op2 A c b ) i ) f ) h ) i' ~> @op2 A ( @op2 A (
  @op2 A ( @op2 A ( @op2 A c h ) b ) f ) i ) i' ). permute. rewrite l
  in v0.  assert ( e * b * d * i * h * i' ~> g * d * b * f * i * i' )
  as k.  assert ( @op2 A ( @op2 A ( @op2 A ( @op2 A ( @op2 A e b ) d )
  i ) h ) i' ~> @op2 A ( @op2 A ( @op2 A ( @op2 A ( @op2 A e h ) i' )
  i ) b ) d ) as k0. permute.  change ( @op2 A ( @op2 A ( @op2 A (
  @op2 A ( @op2 A e b ) d ) i ) h ) i' ~> @op2 A ( @op2 A ( @op2 A (
  @op2 A ( @op2 A g d ) b ) f ) i ) i' ). rewrite k0.  assert ( @op2 A
  ( @op2 A e h ) i' ~> @op2 A ( @op2 A g f ) i' ) as j''. assumption.
  rewrite j''. permute. rewrite k in v0. assumption.
  
  intro u. apply p. intro p'. apply q. intro q'.  destruct p' as [ p'
  j ]. destruct p' as [ i p' ]. destruct q' as [ q' j' ]. destruct q'
  as [ i' q' ].  simpl in *.

  assert ( c * h * b * f * i * i' # g * d * b * f * i * i' ) as v.
  apply azerormultcomp. apply q'. apply azerormultcomp.  apply p'.
  apply azerormultcomp. apply f'. apply azerormultcomp. apply
  b'. assumption.  apply ( pr2 ( acommrng_amult A ) d ). apply ( pr2 (
  acommrng_amult A ) h ).  apply ( pr2 ( acommrng_amult A ) i ). apply
  ( pr2 ( acommrng_amult A ) i' ).

  assert ( c * h * b * f * i * i' ~> a * f * d * h * i * i' ) as k.
  assert ( c * h * b * f * i * i' ~> c * b * i * f * h * i' ) as k0.
  change ( @op2 A ( @op2 A ( @op2 A ( @op2 A ( @op2 A c h ) b ) f ) i
  ) i' ~> @op2 A ( @op2 A ( @op2 A ( @op2 A ( @op2 A c b ) i ) f ) h )
  i' ). permute. rewrite k0. rewrite <- j.  change ( @op2 A ( @op2 A (
  @op2 A ( @op2 A ( @op2 A a d ) i ) f ) h ) i' ~> @op2 A ( @op2 A (
  @op2 A ( @op2 A ( @op2 A a f ) d ) h ) i ) i' ). permute. rewrite k
  in v.  assert ( g * d * b * f * i * i' ~> e * b * d * h * i * i' )
  as l.  assert ( g * d * b * f * i * i' ~> g * f * i' * d * i * b )
  as l0.  change ( @op2 A ( @op2 A ( @op2 A ( @op2 A ( @op2 A g d ) b
  ) f ) i ) i' ~> @op2 A ( @op2 A ( @op2 A ( @op2 A ( @op2 A g f ) i'
  ) d ) i ) b ). permute. rewrite l0. rewrite <- j'.  change (@op2 A (
  @op2 A ( @op2 A ( @op2 A ( @op2 A e h ) i' ) d ) i ) b ~> @op2 A (
  @op2 A ( @op2 A ( @op2 A ( @op2 A e b ) d ) h ) i ) i'
  ). permute. rewrite l in v. assumption.  Defined.

(** We now arrive at the apartness relation on the field of fractions
itself.*)

Definition afldfracapartrel := quotrel afldfracapartiscomprel.

Lemma isirreflafldfracapartrelpre : isirrefl afldfracapartrelpre.
Proof.  intros ab. apply acommrng_airrefl.  Defined.

Lemma issymmafldfracapartrelpre : issymm afldfracapartrelpre.  Proof.
intros ab cd.  apply ( acommrng_asymm A ).  Defined.

Lemma iscotransafldfracapartrelpre : iscotrans afldfracapartrelpre.
Proof.  intros ab cd ef p.  destruct ab as [ a b ]. destruct b as [ b
b' ]. destruct cd as [ c d ] . destruct d as [ d d' ] . destruct ef as
[ e f ]. destruct f as [ f f' ].  assert ( a * f * d # e * b * d ) as
v.  apply azerormultcomp. assumption. assumption.  apply ( (
acommrng_acotrans A ( a * f * d ) ( c * b * f ) ( e * b * d ) ) v
). intro u.  intros P k. apply k.  unfold afldfracapartrelpre in
*. simpl in *. destruct u as [ left | right ].  apply ii1. apply ( pr2
( acommrng_amult A ) f ).  assert ( @op2 A ( @op2 A a f ) d ~> @op2 A
( @op2 A a d ) f ) as i. permute.  change ( @op2 A ( @op2 A a d ) f #
@op2 A ( @op2 A c b ) f ). rewrite <- i. assumption.  apply ii2. apply
( pr2 ( acommrng_amult A ) b ).  assert ( @op2 A ( @op2 A c f ) b ~>
@op2 A ( @op2 A c b ) f ) as i. permute.  change ( @op2 A ( @op2 A c f
) b # @op2 A ( @op2 A e d ) b ). rewrite i.  assert ( @op2 A ( @op2 A
e d ) b ~> @op2 A ( @op2 A e b ) d ) as j. permute.  change ( @op2 A (
@op2 A c b ) f # @op2 A ( @op2 A e d ) b ). rewrite j.  assumption.
Defined.

Lemma isapartafldfracapartrel : isapart afldfracapartrel.  Proof.
intros. split. apply isirreflquotrel. exact (
isirreflafldfracapartrelpre ).  split. apply issymmquotrel. exact (
issymmafldfracapartrelpre ).  apply iscotransquotrel. exact (
iscotransafldfracapartrelpre ).  Defined.

Definition afldfracapart : apart ( commrngfrac A
(aintdomazerosubmonoid A)).  Proof.  intros. unfold apart. split with
afldfracapartrel. exact isapartafldfracapartrel.  Defined.

Lemma isbinapartlafldfracop1 : isbinopapartl afldfracapart op1.
Proof.  intros. unfold isbinopapartl.  assert ( forall a b c :
commrngfrac A ( aintdomazerosubmonoid A ), isaprop ( pr1
(afldfracapart) ( commrngfracop1 A ( aintdomazerosubmonoid A ) a b) (
commrngfracop1 A ( aintdomazerosubmonoid A ) a c) -> pr1
(afldfracapart ) b c) ) as int.  intros a b c. apply impred. intro
p. apply ( pr1 ( afldfracapart ) b c ).  apply ( setquotuniv3prop _ (
fun a b c => hProppair _ ( int a b c ) ) ).  intros ab cd ef
p. destruct ab as [ a b ]. destruct b as [ b b' ].  destruct cd as [ c
d ]. destruct d as [ d d' ]. destruct ef as [ e f ].  destruct f as [
f f' ]. unfold afldfracapart in *. simpl.  unfold
afldfracapartrel. unfold quotrel. rewrite setquotuniv2comm. unfold
afldfracapartrelpre.  simpl.

  assert ( afldfracapartrelpre ( dirprodpair ( @op1 A ( @op2 A d a ) (
  @op2 A b c ) ) ( @op ( aintdomazerosubmonoid A ) ( tpair b b' ) (
  tpair d d' ) ) ) ( dirprodpair ( @op1 A ( @op2 A f a ) ( @op2 A b e
  ) ) ( @op ( aintdomazerosubmonoid A ) ( tpair b b' ) ( tpair f f' )
  ) ) ) as u. apply p.  unfold afldfracapartrelpre in u. simpl in
  u. rewrite 2! ( @rngrdistr A ) in u.  repeat rewrite <- rngassoc2 in
  u. assert ( (@op2 (pr1rng (commrngtorng (acommrngtocommrng
  (pr1aintdom A)))) (@op2 (pr1rng (commrngtorng (acommrngtocommrng
  (pr1aintdom A)))) (@op2 (@pr1 setwith2binop (fun X : setwith2binop
  => @iscommrngops (pr1setwith2binop X) (@op1 X) (@op2 X))
  (acommrngtocommrng (pr1aintdom A))) d a) b) f) ~> (@op2 (pr1rng
  (commrngtorng (acommrngtocommrng (pr1aintdom A)))) (@op2 (pr1rng
  (commrngtorng (acommrngtocommrng (pr1aintdom A)))) (@op2 (@pr1
  setwith2binop (fun X : setwith2binop => @iscommrngops
  (pr1setwith2binop X) (@op1 X) (@op2 X)) (acommrngtocommrng
  (pr1aintdom A))) f a) b) d) ) as i.  permute. rewrite i in u. assert
  ( (@op2 (pr1rng (commrngtorng (acommrngtocommrng (pr1aintdom A))))
  (@op2 (pr1rng (commrngtorng (acommrngtocommrng (pr1aintdom A))))
  (@op2 (@pr1 setwith2binop (fun X : setwith2binop => @iscommrngops
  (pr1setwith2binop X) (@op1 X) (@op2 X)) (acommrngtocommrng
  (pr1aintdom A))) b c) b) f) ~> (@op2 (pr1rng (commrngtorng
  (acommrngtocommrng (pr1aintdom A)))) (@op2 (pr1rng (commrngtorng
  (acommrngtocommrng (pr1aintdom A)))) (@op2 (@pr1 setwith2binop (fun
  X : setwith2binop => @iscommrngops (pr1setwith2binop X) (@op1 X)
  (@op2 X)) (acommrngtocommrng (pr1aintdom A))) c f) b) b) ) as
  j. permute. rewrite j in u.  assert ( (@op2 (pr1rng (commrngtorng
  (acommrngtocommrng (pr1aintdom A)))) (@op2 (pr1rng (commrngtorng
  (acommrngtocommrng (pr1aintdom A)))) (@op2 (@pr1 setwith2binop (fun
  X : setwith2binop => @iscommrngops (pr1setwith2binop X) (@op1 X)
  (@op2 X)) (acommrngtocommrng (pr1aintdom A))) b e) b) d) ~> (@op2
  (pr1rng (commrngtorng (acommrngtocommrng (pr1aintdom A)))) (@op2
  (pr1rng (commrngtorng (acommrngtocommrng (pr1aintdom A)))) (@op2
  (@pr1 setwith2binop (fun X : setwith2binop => @iscommrngops
  (pr1setwith2binop X) (@op1 X) (@op2 X)) (acommrngtocommrng
  (pr1aintdom A))) e d) b) b) ) as j'. permute. rewrite j' in u.
  apply ( pr2 ( acommrng_amult A ) b ).  apply ( pr2 ( acommrng_amult
  A ) b ).  apply ( pr1 ( acommrng_aadd A) ( f * a * b * d )
  ). assumption.  Defined.

Lemma isbinapartrafldfracop1 : isbinopapartr afldfracapart op1.
Proof.  intros a b c. rewrite ( rngcomm1 ). rewrite ( rngcomm1 _ c
). apply isbinapartlafldfracop1.  Defined.

Lemma isbinapartlafldfracop2 : isbinopapartl afldfracapart op2.
Proof.  intros. unfold isbinopapartl.  assert ( forall a b c :
commrngfrac A ( aintdomazerosubmonoid A ), isaprop ( pr1
(afldfracapart ) ( commrngfracop2 A ( aintdomazerosubmonoid A ) a b) (
commrngfracop2 A ( aintdomazerosubmonoid A ) a c) -> pr1
(afldfracapart ) b c) ) as int.  intros a b c. apply impred. intro
p. apply ( pr1 ( afldfracapart ) b c ).  apply ( setquotuniv3prop _ (
fun a b c => hProppair _ ( int a b c ) ) ).  intros ab cd ef p.
  
  destruct ab as [ a b ]. destruct b as [ b b' ].  destruct cd as [ c
  d ]. destruct d as [ d d' ]. destruct ef as [ e f ].  destruct f as
  [ f f' ].
  
  assert ( afldfracapartrelpre ( dirprodpair ( ( a * c ) ) ( @op (
  aintdomazerosubmonoid A ) ( tpair b b' ) ( tpair d d' ) ) ) (
  dirprodpair ( a * e ) ( @op ( aintdomazerosubmonoid A ) ( tpair b b'
  ) ( tpair f f' ) ) ) ) as u. apply p.  unfold afldfracapart in
  *. simpl. unfold afldfracapartrel.  unfold quotrel. rewrite (
  setquotuniv2comm ( eqrelcommrngfrac A ( aintdomazerosubmonoid A ) )
  ).  unfold afldfracapartrelpre in *. simpl. simpl in u.  apply ( pr2
  ( acommrng_amult A ) a ). apply ( pr2 ( acommrng_amult A ) b ).


  assert ( c * f * a * b ~> (@op2 (@pr1 setwith2binop (fun X :
  setwith2binop => @iscommrngops (pr1setwith2binop X) (@op1 X) (@op2
  X)) (acommrngtocommrng (pr1aintdom A))) (@op2 (@pr1 setwith2binop
  (fun X : setwith2binop => @iscommrngops (pr1setwith2binop X) (@op1
  X) (@op2 X)) (acommrngtocommrng (pr1aintdom A))) a c) (@op2 (@pr1
  setwith2binop (fun X : setwith2binop => @iscommrngops
  (pr1setwith2binop X) (@op1 X) (@op2 X)) (acommrngtocommrng
  (pr1aintdom A))) b f)) ) as i. change ( c * f * a * b ~> a * c * ( b
  * f ) ). permute.  change ( c * f * a * b # e * d * a * b ). rewrite
  i.  assert ( e * d * a * b ~> (@op2 (@pr1 setwith2binop (fun X :
  setwith2binop => @iscommrngops (pr1setwith2binop X) (@op1 X) (@op2
  X)) (acommrngtocommrng (pr1aintdom A))) (@op2 (@pr1 setwith2binop
  (fun X : setwith2binop => @iscommrngops (pr1setwith2binop X) (@op1
  X) (@op2 X)) (acommrngtocommrng (pr1aintdom A))) a e) (@op2 (@pr1
  setwith2binop (fun X : setwith2binop => @iscommrngops
  (pr1setwith2binop X) (@op1 X) (@op2 X)) (acommrngtocommrng
  (pr1aintdom A))) b d)) ) as i'. change ( e * d * a * b ~> a * e * (
  b * d ) ). permute.  rewrite i'. assumption.  Defined.

Lemma isbinapartrafldfracop2 : isbinopapartr (afldfracapart ) op2.
Proof.  intros a b c. rewrite rngcomm2. rewrite ( rngcomm2 _ c
). apply isbinapartlafldfracop2.  Defined.

Definition afldfrac0 : acommrng.  Proof.  intros. split with (
commrngfrac A ( aintdomazerosubmonoid A ) ).  split with (
afldfracapart ). split.  split. apply ( isbinapartlafldfracop1
). apply ( isbinapartrafldfracop1 ).  split. apply (
isbinapartlafldfracop2 ). apply ( isbinapartrafldfracop2 ).  Defined.

Definition afldfracmultinvint ( ab : dirprod A ( aintdomazerosubmonoid
A ) ) ( is : afldfracapartrelpre ab ( dirprodpair ( @rngunel1 A ) (
unel ( aintdomazerosubmonoid A ) )) ) : dirprod A (
aintdomazerosubmonoid A ).  Proof.  intros. destruct ab as [ a b
]. destruct b as [ b b' ].  split with b. simpl in is. split with
a. unfold afldfracapartrelpre in is. simpl in is.  change ( a # 0
). rewrite ( @rngmult0x A ) in is. rewrite ( @rngrunax2 A ) in
is. assumption.  Defined.

Definition afldfracmultinv ( a : afldfrac0 ) ( is : a # 0 ) :
multinvpair afldfrac0 a.  Proof.  intros. assert ( forall b :
afldfrac0, isaprop ( b # 0 -> multinvpair afldfrac0 b ) ) as int.
intros. apply impred. intro p. apply ( isapropmultinvpair afldfrac0 ).
assert ( forall b : afldfrac0, b # 0 -> multinvpair afldfrac0 b ) as
p.  apply ( setquotunivprop _ ( fun x0 => hProppair _ ( int x0 ) ) ).
intros bc q.  destruct bc as [ b c ].  assert ( afldfracapartrelpre (
dirprodpair b c ) ( dirprodpair ( @rngunel1 A ) ( unel (
aintdomazerosubmonoid A ) ) ) ) as is'.  apply q. split with
(setquotpr (eqrelcommrngfrac A (aintdomazerosubmonoid A)) (
afldfracmultinvint ( dirprodpair b c ) is' ) ).
  
  split.  change ( setquotpr ( eqrelcommrngfrac A (
  aintdomazerosubmonoid A ) ) ( dirprodpair ( @op2 A ( pr1 (
  afldfracmultinvint ( dirprodpair b c ) is' ) ) b ) ( @op (
  aintdomazerosubmonoid A ) ( pr2 ( afldfracmultinvint ( dirprodpair b
  c ) is' ) ) c ) ) ~> ( commrngfracunel2 A ( aintdomazerosubmonoid A
  ) ) ).  apply iscompsetquotpr. unfold commrngfracunel2int.  destruct
  c as [ c c' ]. simpl.  apply total2tohexists. split with (
  carrierpair ( fun x : pr1 A => x # 0 ) 1 ( pr1 ( pr2 A ) ) ).
  simpl. rewrite 3! ( @rngrunax2 A ).  rewrite ( @rnglunax2 A ). apply
  ( @rngcomm2 A ).

  change ( setquotpr ( eqrelcommrngfrac A ( aintdomazerosubmonoid A )
 ) ( dirprodpair ( @op2 A b ( pr1 ( afldfracmultinvint ( dirprodpair b
 c ) is' ) ) ) ( @op ( aintdomazerosubmonoid A ) c ( pr2 (
 afldfracmultinvint ( dirprodpair b c ) is' ) ) ) ) ~> (
 commrngfracunel2 A ( aintdomazerosubmonoid A ) ) ).  apply
 iscompsetquotpr. destruct c as [ c c' ]. simpl. apply
 total2tohexists.  split with ( carrierpair ( fun x : pr1 A => x # 0 )
 1 ( pr1 ( pr2 A ) ) ).  simpl. rewrite 3! ( @rngrunax2 A ). rewrite (
 @rnglunax2 A ). apply ( @rngcomm2 A ).  apply p. assumption.
 Defined.
 
Theorem afldfracisafld : isaafield afldfrac0.  Proof.  intros. split.
change ( ( afldfracapartrel ) ( @rngunel2 ( commrngfrac A (
aintdomazerosubmonoid A ) ) ) ( @rngunel1 ( commrngfrac A (
aintdomazerosubmonoid A ) ) ) ).  unfold afldfracapartrel. cut ( (
@op2 A ( @rngunel2 A ) ( @rngunel2 A ) ) # ( @op2 A ( @rngunel1 A ) (
@rngunel2 A ) ) ).  intro v. apply v. rewrite 2! ( @rngrunax2 A
). apply A.
  
  intros a p. apply afldfracmultinv. assumption.  Defined.
 
Definition afldfrac := afldpair afldfrac0 afldfracisafld.

End aint.

Close Scope rng_scope.
(** END OF FILE *)
\end{verbatim}
}
}\end{multicols}

\subsection{The file \texttt{zmodp.v}}

\begin{multicols}{2}{
{
\tiny

}
}\end{multicols}

\subsection{The file \texttt{padics.v}}

\begin{multicols}{2}{
{
\tiny
%
}
}\end{multicols}

\end{landscape}

{\small
\begin{multicols}{3}{
    \begin{center}
      \'Alvaro Pelayo \\
      School of Mathematics\\
      Institute for Advanced Study\\
      Einstein Drive, Princeton\\
      NJ 08540 USA.\\
      \vspace{.2cm}
      and\\
      \vspace{.2cm}
      \noindent
      Washington University\\  
      Mathematics Department \\
      One Brookings Drive, Campus Box 1146\\
      St Louis\\
      MO 63130-4899, USA.\\
      {\em E-mail}: \texttt{apelayo@math.wustl.edu}
    \end{center}
\columnbreak
\begin{center}
  Vladimir Voevodsky\\
  School of Mathematics \\
  Institute for Advanced Study\\
  Einstein Drive, Princeton\\
  NJ 08540 USA. \\
  {\em E-mail}: \texttt{vladimir@math.ias.edu}
\end{center}
\columnbreak
\begin{center}
  Michael A. Warren\\
  School of Mathematics\\
  Institute for Advanced Study\\
  Einstein Drive, Princeton\\
  NJ 08540 USA.\\
  {\em E-mail}: \texttt{mwarren@math.ias.edu}
\end{center}
}\end{multicols}
}
\end{document}